\documentclass[twoside,leqno]{article}
\usepackage{amsmath}
\usepackage{amsthm}
\usepackage{amssymb}
\usepackage{amscd}
\usepackage{graphicx}
\usepackage{epsfig}

\usepackage{latexsym}
\usepackage{amsfonts}
\input xy
\xyoption{all} \tolerance=500

-\headheight\advance\topmargin by -\headsep

\numberwithin{equation}{section} \allowdisplaybreaks

\newtheorem{theorem}{Theorem}[section]

\theoremstyle{definition}
\newtheorem{definition}{Definition}[section]
\newtheorem{example}{Example} [section]

\newtheorem{theo}{Theorem}
\begin{document}
\font\black=cmbx10 \font\sblack=cmbx7 \font\ssblack=cmbx5
\font\blackital=cmmib10  \skewchar\blackital='177
\font\sblackital=cmmib7 \skewchar\sblackital='177
\font\ssblackital=cmmib5 \skewchar\ssblackital='177
\font\sanss=cmss10 \font\ssanss=cmss8 
\font\sssanss=cmss8 scaled 600 \font\blackboard=msbm10
\font\sblackboard=msbm7 \font\ssblackboard=msbm5
\font\caligr=eusm10 \font\scaligr=eusm7 \font\sscaligr=eusm5
\font\blackcal=eusb10 \font\fraktur=eufm10 \font\sfraktur=eufm7
\font\ssfraktur=eufm5 \font\blackfrak=eufb10

\font\bsymb=cmsy10 scaled\magstep2
\def\all#1{\setbox0=\hbox{\lower1.5pt\hbox{\bsymb
       \char"38}}\setbox1=\hbox{$_{#1}$} \box0\lower2pt\box1\;}
\def\exi#1{\setbox0=\hbox{\lower1.5pt\hbox{\bsymb \char"39}}
       \setbox1=\hbox{$_{#1}$} \box0\lower2pt\box1\;}

\def\mi#1{{\fam1\relax#1}}
\def\tx#1{{\fam0\relax#1}}

\newfam\bifam
\textfont\bifam=\blackital \scriptfont\bifam=\sblackital
\scriptscriptfont\bifam=\ssblackital
\def\bi#1{{\fam\bifam\relax#1}}

\newfam\blfam
\textfont\blfam=\black \scriptfont\blfam=\sblack
\scriptscriptfont\blfam=\ssblack
\def\rbl#1{{\fam\blfam\relax#1}}

\newfam\bbfam
\textfont\bbfam=\blackboard \scriptfont\bbfam=\sblackboard
\scriptscriptfont\bbfam=\ssblackboard
\def\bb#1{{\fam\bbfam\relax#1}}

\newfam\ssfam
\textfont\ssfam=\sanss \scriptfont\ssfam=\ssanss
\scriptscriptfont\ssfam=\sssanss
\def\sss#1{{\fam\ssfam\relax#1}}

\newfam\clfam
\textfont\clfam=\caligr \scriptfont\clfam=\scaligr
\scriptscriptfont\clfam=\sscaligr
\def\cl#1{{\fam\clfam\relax#1}}

\newfam\frfam
\textfont\frfam=\fraktur \scriptfont\frfam=\sfraktur
\scriptscriptfont\frfam=\ssfraktur
\def\fr#1{{\fam\frfam\relax#1}}

\def\cb#1{\hbox{$\fam\gpfam\relax#1\textfont\gpfam=\blackcal$}}

\def\hpb#1{\setbox0=\hbox{${#1}$}
    \copy0 \kern-\wd0 \kern.2pt \box0}
\def\vpb#1{\setbox0=\hbox{${#1}$}
    \copy0 \kern-\wd0 \raise.08pt \box0}

\def\pmb#1{\setbox0\hbox{${#1}$} \copy0 \kern-\wd0 \kern.2pt \box0}
\def\pmbb#1{\setbox0\hbox{${#1}$} \copy0 \kern-\wd0
      \kern.2pt \copy0 \kern-\wd0 \kern.2pt \box0}
\def\pmbbb#1{\setbox0\hbox{${#1}$} \copy0 \kern-\wd0
      \kern.2pt \copy0 \kern-\wd0 \kern.2pt
    \copy0 \kern-\wd0 \kern.2pt \box0}
\def\pmxb#1{\setbox0\hbox{${#1}$} \copy0 \kern-\wd0
      \kern.2pt \copy0 \kern-\wd0 \kern.2pt
      \copy0 \kern-\wd0 \kern.2pt \copy0 \kern-\wd0 \kern.2pt \box0}
\def\pmxbb#1{\setbox0\hbox{${#1}$} \copy0 \kern-\wd0 \kern.2pt
      \copy0 \kern-\wd0 \kern.2pt
      \copy0 \kern-\wd0 \kern.2pt \copy0 \kern-\wd0 \kern.2pt
      \copy0 \kern-\wd0 \kern.2pt \box0}

\def\cdotss{\mathinner{\cdotp\cdotp\cdotp\cdotp\cdotp\cdotp\cdotp
        \cdotp\cdotp\cdotp\cdotp\cdotp\cdotp\cdotp\cdotp\cdotp\cdotp
        \cdotp\cdotp\cdotp\cdotp\cdotp\cdotp\cdotp\cdotp\cdotp\cdotp
        \cdotp\cdotp\cdotp\cdotp\cdotp\cdotp\cdotp\cdotp\cdotp\cdotp}}

\font\frak=eufm10 scaled\magstep1 \font\fak=eufm10 scaled\magstep2
\font\fk=eufm10 scaled\magstep3 \font\scriptfrak=eufm10
\font\tenfrak=eufm10


\mathchardef\za="710B  
\mathchardef\zb="710C  
\mathchardef\zg="710D  
\mathchardef\zd="710E  
\mathchardef\zve="710F 
\mathchardef\zz="7110  
\mathchardef\zh="7111  
\mathchardef\zvy="7112 
\mathchardef\zi="7113  
\mathchardef\zk="7114  
\mathchardef\zl="7115  
\mathchardef\zm="7116  
\mathchardef\zn="7117  
\mathchardef\zx="7118  
\mathchardef\zp="7119  
\mathchardef\zr="711A  
\mathchardef\zs="711B  
\mathchardef\zt="711C  
\mathchardef\zu="711D  
\mathchardef\zvf="711E 
\mathchardef\zq="711F  
\mathchardef\zc="7120  
\mathchardef\zw="7121  
\mathchardef\ze="7122  
\mathchardef\zy="7123  
\mathchardef\zf="7124  
\mathchardef\zvr="7125 
\mathchardef\zvs="7126 
\mathchardef\zf="7127  
\mathchardef\zG="7000  
\mathchardef\zD="7001  
\mathchardef\zY="7002  
\mathchardef\zL="7003  
\mathchardef\zX="7004  
\mathchardef\zP="7005  
\mathchardef\zS="7006  
\mathchardef\zU="7007  
\mathchardef\zF="7008  
\mathchardef\zW="700A  

\newcommand{\be}{\begin{equation}}
\newcommand{\ee}{\end{equation}}
\newcommand{\ra}{\rightarrow}
\newcommand{\lra}{\longrightarrow}
\newcommand{\bea}{\begin{eqnarray}}
\newcommand{\eea}{\end{eqnarray}}
\newcommand{\beas}{\begin{eqnarray*}}
\newcommand{\eeas}{\end{eqnarray*}}
\def\*{{\textstyle *}}
\newcommand{\R}{{\mathbb R}}
\newcommand{\T}{{\mathbb T}}
\newcommand{\C}{{\mathbb C}}
\newcommand{\unit}{{\mathbf 1}}
\newcommand{\SL}{SL(2,\C)}
\newcommand{\Sl}{sl(2,\C)}
\newcommand{\SU}{SU(2)}
\newcommand{\su}{su(2)}
\def\ssT{\sss T}
\newcommand{\G}{{\goth g}}
\newcommand{\D}{{\rm d}}
\newcommand{\Df}{{\rm d}^\zF}
\newcommand{\de}{\,{\stackrel{\rm def}{=}}\,}
\newcommand{\we}{\wedge}
\newcommand{\nn}{\nonumber}
\newcommand{\ot}{\otimes}
\newcommand{\s}{{\textstyle *}}
\newcommand{\ts}{T^\s}
\newcommand{\oX}{\stackrel{o}{X}}
\newcommand{\oD}{\stackrel{o}{D}}
\newcommand{\obD}{\stackrel{o}{\bD}}
\newcommand{\pa}{\partial}
\newcommand{\ti}{\times}
\newcommand{\A}{{\cal A}}
\newcommand{\Li}{{\cal L}}
\newcommand{\ka}{\mathbb{K}}
\newcommand{\find}{\mid}
\newcommand{\ad}{{\rm ad}}
\newcommand{\rS}{]^{SN}}
\newcommand{\rb}{\}_P}
\newcommand{\p}{{\sf P}}
\newcommand{\h}{{\sf H}}
\newcommand{\X}{{\cal X}}
\newcommand{\I}{\,{\rm i}\,}
\newcommand{\rB}{]_P}
\newcommand{\Ll}{{\pounds}}
\def\lna{\lbrack\! \lbrack}
\def\rna{\rbrack\! \rbrack}
\def\rnaf{\rbrack\! \rbrack_\zF}
\def\rnah{\rbrack\! \rbrack\,\hat{}}
\def\lbo{{\lbrack\!\!\lbrack}}
\def\rbo{{\rbrack\!\!\rbrack}}
\def\lan{\langle}
\def\ran{\rangle}
\def\zT{{\cal T}}
\def\tU{\tilde U}
\def\ati{{\stackrel{a}{\times}}}
\def\sti{{\stackrel{sv}{\times}}}
\def\aot{{\stackrel{a}{\ot}}}
\def\sati{{\stackrel{sa}{\times}}}
\def\saop{{\stackrel{sa}{\op}}}
\def\bwa{{\stackrel{a}{\bigwedge}}}
\def\svop{{\stackrel{sv}{\oplus}}}
\def\saot{{\stackrel{sa}{\otimes}}}
\def\cti{{\stackrel{cv}{\times}}}
\def\cop{{\stackrel{cv}{\oplus}}}
\def\dra{{\stackrel{\xd}{\ra}}}
\def\bdra{{\stackrel{\bd}{\ra}}}
\def\bAff{\mathbf{Aff}}
\def\Aff{\sss{Aff}}
\def\bHom{\mathbf{Hom}}
\def\Hom{\sss{Hom}}
\def\bt{{\boxtimes}}
\def\sot{{\stackrel{sa}{\ot}}}
\def\bp{{\boxplus}}
\def\op{\oplus}
\def\bwak{{\stackrel{a}{\bigwedge}\!{}^k}}
\def\aop{{\stackrel{a}{\oplus}}}
\def\ix{\operatorname{i}}
\def\V{{\cal V}}
\def\cD{{\cal D}}
\def\cC{{\cal C}}
\def\cE{{\cal E}}
\def\cL{{\cal L}}
\def\cN{{\cal N}}
\def\cR{{\cal R}}
\def\cJ{{\cal J}}
\def\cT{{\cal T}}
\def\cH{{\cal H}}
\def\bA{\mathbf{A}}
\def\bI{\mathbf{I}}
\def\wh{\widehat}
\def\wt{\widetilde}
\def\ol{\overline}
\def\ul{\underline}
\def\Sec{\sss{Sec}}
\def\Lin{\sss{Lin}}
\def\ader{\sss{ADer}}
\def\ado{\sss{ADO^1}}
\def\adoo{\sss{ADO^0}}
\def\AS{\sss{AS}}
\def\bAS{\sss{AS}}
\def\bLS{\sss{LS}}
\def\bAP{\sss{AV}}
\def\bLP{\sss{LP}}
\def\AP{\sss{AP}}
\def\LP{\sss{LP}}
\def\LS{\sss{LS}}
\def\Z{\mathbf{Z}}
\def\oZ{\overline{\bZ}}
\def\oA{\overline{\bA}}
\def\cim{{C^\infty(M)}}
\def\de{{\cal D}^1}
\def\la{\langle}
\def\ran{\rangle}
\def\by{{\bi y}}
\def\bs{{\bi s}}
\def\bc{{\bi c}}
\def\bd{{\bi d}}
\def\bh{{\bi h}}
\def\bD{{\bi D}}
\def\bY{{\bi Y}}
\def\bX{{\bi X}}
\def\bL{{\bi L}}
\def\bV{{\bi V}}
\def\bW{{\bi W}}
\def\bS{{\bi S}}
\def\bT{{\bi T}}
\def\bC{{\bi C}}
\def\bE{{\bi E}}
\def\bF{{\bi F}}
\def\bP{{\bi P}}
\def\bp{{\bi p}}
\def\bz{{\bi z}}
\def\bZ{{\bi Z}}
\def\bq{{\bi q}}
\def\bQ{{\bi Q}}
\def\bx{{\bi x}}

\def\sA{{\sss A}}
\def\sC{{\sss C}}
\def\sD{{\sss D}}
\def\sG{{\sss G}}
\def\sH{{\sss H}}
\def\sI{{\sss I}}
\def\sJ{{\sss J}}
\def\sK{{\sss K}}
\def\sL{{\sss L}}
\def\sO{{\sss O}}
\def\sP{{\sss P}}
\def\sPh{{\sss P\sss h}}
\def\sT{{\sss T}}
\def\sV{{\sss V}}
\def\sR{{\sss R}}
\def\sS{{\sss S}}
\def\sE{{\sss E}}
\def\sF{{\sss F}}
\def\st{{\sss t}}
\def\sg{{\sss g}}
\def\sx{{\sss x}}
\def\sv{{\sss v}}
\def\sw{{\sss w}}
\def\sQ{{\sss Q}}
\def\sj{{\sss j}}
\def\sq{{\sss q}}
\def\xa{\tx{a}}
\def\xc{\tx{c}}
\def\xd{\tx{d}}
\def\xi{\tx{i}}
\def\xD{\tx{D}}
\def\xV{\tx{V}}
\def\xF{\tx{F}}


\setcounter{page}{1} \thispagestyle{empty}




\bigskip

\bigskip

\title{AV-differential geometry: Euler-Lagrange equations\thanks{Research
supported by the Polish Ministry of Scientific Research and
Information Technology under the grant No. 2 P03A 036 25.}}

        \author{
        Katarzyna  Grabowska$^1$, Janusz Grabowski$^2$, Pawe\l\ Urba\'nski$^1$\\
        \\
         $^1$ {\it Physics Department}\\
                {\it University of Warsaw} \\
         $^2$ {\it Institute of Mathematics}\\
                {\it Polish Academy of Sciences}
                }
\date{}
\maketitle
\begin{abstract}
A general, consistent and complete framework for geometrical
formulation of mechanical systems is proposed, based on certain
structures on affine bundles (affgebroids) that generalize Lie
algebras and Lie algebroids. This scheme covers and unifies
various geometrical approaches to mechanics in the Lagrangian and
Hamiltonian pictures, including time-dependent lagrangians and
hamiltonians. In our approach, lagrangians and hamiltonians are,
in general, sections of certain $\R$-principal bundles, and the
solutions of analogs of Euler-Lagrange equations are curves in
certain affine bundles. The correct geometrical and
frame-independent description of Newtonian Mechanics is of this
type.


\end{abstract}
\section{Introduction} In earlier papers \cite{GGU2,Ur1} we
developed a geometrical theory (AV-geometry) in which functions on
a mani\-fold are replaced by sections of $\R$-principal bundles
(AV-bundles). It was used for frame-independent formulation of a
number of physical theories. The general geometrical concepts and
tools described in \cite{GGU2} were then applied to obtain the
Hamiltonian picture in this affine setting.

In the present paper, in turn, we are concentrated on the
Lagrangian picture and we derive analogs of the Euler-Lagrange
equations in the AV-bundle setting. This is done completely
geometrically and intrinsically, however with no reference to any
variational calculus (which we plan to to study in our forthcoming
paper). We use the framework of what we call {\it special
affgebroids} which is more general than that of {\it Lie
affgebroids} used recently by other authors \cite{IMPS}. We get a
larger class of Euler-Lagrange equations and larger class of
possible models for physical theories which cover the ones
considered in \cite{IMPS} as particular examples. The special
affgebroids are related to Lie affgebroids in the way the general
{\it algebroids} used in \cite{GGU3} are related to Lie
algebroids. The other difference is that we do not use
prolongations of Lie affgebroids (nor affgebroids) that
simplifies, in our opinion, the whole picture. This is possible,
since in our method we do not follow the Klein's ideas \cite{Kl}
for geometric construction of Euler-Lagrange equations but rather
some ideas due to W.~M.~Tulczyjew \cite{Tu2}. The dynamics
obtained in this way is therefore implicit but we think it is the
nature of the problem and this approach has the advantage that
regularity of the Lagrangian plays no role in the main
construction.

Note that the idea of using affine bundles for geometrically
correct description of mechanical systems is not new and goes back
to \cite{TU1,Tu,UR}. A similar approach to time-dependent
non-relativistic mechanics in the Lagrange formulation has been
recently developed by other authors \cite{MMS0,MMS1,MPL,MVB,MMS}.

The paper is organized as follows. In Section 2 we recall
rudiments of the AV-geometry that will be needed in the sequel.
Section 3 is devoted to the presentation of the concept of double
affine bundle which is the fundamental concept in our approach. In
Section 4 we recall the Lagrangian and Hamiltonian formalisms for
general algebroids developed in \cite{GGU3}. On this fundamentals,
special affgebroids as certain morphisms of double affine bundles
and the corresponding Lagrangian and Hamiltonian formalisms are
constructed and studied in Section 5. In particular, we derive
analogs of Euler-Lagrange equations. We end up with some examples
in Section 6.

\section{Rudiments of the AV-geometry}
We refer to \cite{GGU2} (see also \cite{GGU1,Ur1}) for a
development of the AV-geometry. To make the paper self-contained,
we recall basic notions and facts here.

\subsection{Affine spaces and affine bundles}
An {\it affine space} is a triple $(A,V,\za)$, where $A$ is a set,
$V$ is a vector space over a field $\ka$ and $\za$ is a mapping
$\za \colon A \times A\rightarrow V$ such that
      \begin{itemize}
   \item $\za(a_3,a_2) + \za(a_2,a_1) + \za(a_1,a_3) = 0$;
   \item the mapping $\za(\cdot,a) \colon A \rightarrow V$ is bijective for
each $a \in A$.
      \end{itemize}
One writes usually the map $\za$ as a difference $\za(a,a')=a-a'$.
Equivalently, one can add vectors from $V$ to points of $A$:
$\za(\cdot,a)^{-1}(v)=a+v$. We shall also write simply $A$ to
denote the affine space $(A,V,\za)$ and $\sV(A)$ to denote its
{\it model vector space} $V$. The {\it dimension of an affine
space} $A$ we will call the dimension of its model vector space
$\sV(A)$. A mapping $\zf$ from the affine space $A_1$ to the
affine space $A_2$ is {\it affine} if there exists a linear map
$\zf_\sv: \sV(A_1)\ra\sV(A_2)$ such that, for every $a\in A_1$ and
$u\in \sV(A_1)$,
    $$\zf(a+u)=\zf(a)+\zf_\sV(u).$$
The mapping $\zf_\sv$ is called the {\it linear part} of $\zf$. We
will need also {\it multi-affine}, especially {\it bi-affine},
maps. Let $A$, $A_1$, $A_2$ be affine spaces. A map
    $$\Phi: A_1\times A_2 \longrightarrow A$$
is called {\it bi-affine} if it is affine with respect to every
argument separately. One can also define linear parts of a
bi-affine map. By $\Phi_{\sv}^1$ we will denote the linear-affine
map
    $$\Phi_{\sv}^1: \sV(A_1)\times A_2\longrightarrow \sV(A),$$
where the linear part is taken with respect to the first argument.
Similarly, for the second argument we have the affine-linear map
    $$\Phi_{\sv}^2: A_1\times \sV(A_2)\longrightarrow \sV(A)$$
and, finally, the {\it bilinear part of\ } $\Phi$:
    $$\Phi_{\sv}:\sV(A_1)\times \sV(A_2)\longrightarrow \sV(A).$$
 For an affine space
$A$ we define its {\it vector dual} $A^\dag$ as the set of all
affine functions  $\phi:A\ra\ka$. The dimension of $A^\dag$ is
greater by $1$ than the dimension of $A$. By $1_A$ we will denote
the element of $A^\dag$ being the constant function on $A$ equal
to $1$. Observe that $A$ is naturally embedded in the space
$\widehat A=(A^\dag)^\ast$. The affine space $A$ can be then
identified with the one-codimensional affine subspace of $\widehat
A$ of those linear functions on $A^\dag$ that evaluated on $1_A$
give $1$. Similarly, $\sV(A)$ can be identified with the vector
subspace of $\widehat A$ of those elements that evaluated on $1_A$
give $0$. The above observation justifies the name {\it vector
hull of $A$}, for $\widehat A$.

In the following we will widely use {\it affine bundles} which are
smooth, locally trivial bundles of affine spaces. Here $\ka=\R$
and the notation will be, in principle, the same for affine spaces
and affine bundles. For instance, $\sV(A)$ denotes the vector
bundle which is the model for an affine bundle $\zz:A\ra M$ over a
base manifold $M$. By $\Sec$ we denote the spaces of sections,
e.g. $\Sec(A)$ is the affine space of sections of the affine
bundle $\zz:A\ra M$. The difference of two sections $a$, $a'$ of
the bundle $A$ is a section of $\sV(A)$. Equivalently, we can add
sections of $V(A)$ to sections of $A$:
$$\Sec(A)\ti\Sec(\sV(A))\ni(a,v)\mapsto a+v\in\Sec(A).$$
Every $a\in \Sec(A)$ induces a `linearization'
$$I_a:A\ra\sV(A),\quad a'_p\mapsto a'_p-a(p).$$
The {\it bundle of affine morphisms } (over the identity on the
base) $\Aff_M(A_1;A_2)$ from $A_1$ to $A_2$ is a bundle which
fibers are spaces of affine maps from fibers of $A_1$ to fibers of
$A_2$ over the same point in $M$. The space of sections of this
bundle we will denote by $\Aff(A_1;A_2)$. It consists of actual
morphisms of affine bundles over the identity on the base.
Morphisms of affine bundles over an arbitrary map of the bases are
defined in the obvious way, analogous to that for vector bundles.
The bundle $\Aff_M(A_1;A_2)$ is an affine bundle modelled on
$\Aff_M(A_1;\sV(A_2))$ and $\Aff(A_1;A_2)$ is an affine space
modelled on $\Aff(A_1;\sV(A_2))$. Like for affine maps, we can
define a {\it linear part of an affine morphism}: if $\zf\in
\Aff(A_1;A_2)$, then $\zf_\sv\in \Hom(\sV(A_1);\sV(A_2))$ and for
any $a\in A_1$ and $u\in V(A_1)$, both over the same point in $M$,
we have
    $$\zf_\sv(u)=\zf(a+u)-\zf(a).$$
For an affine bundle $A$ we define also its {\it vector dual}
$A^\dag=\Aff_M(A;M\ti\R)$ and its {\it vector hull}, i.e. the
vector bundle $\wh{A}=(A^\dag)^\ast$.

\subsection{Special affine bundles} A vector space with a
distinguished non-zero vector will be called {\it a special vector
space}. A canonical example of a special vector space is $(\R,1)$.
Another example of a special vector space is the vector dual
$\bA^\dag=(A^\dag, 1_A)$ of an affine space $A$. A {\it special
affine} space is an affine space modelled on a special vector
space. Similarly, a vector bundle $\bV$ with a distinguished
nowhere-vanishing section $v$ will be called {\it a special vector
bundle.} An affine bundle modelled on a special vector bundle will
be called {\it a special affine bundle}. If $\bA$ is a special
affine bundle, then $\wh{\bA}$ is canonically a special vector
bundle. If $\bI$ denotes the special vector space $(\R,1)$
understood as a special affine space, then $M\times \bI$ is a
special affine bundle over $M$. With some abuse of notation this
bundle will be also denoted by $\bI$. If $\bA_i$ denotes a special
affine bundle over $M$ modelled on a special vector bundle
$\bV_i=\sV(\bA_i)$ with the distinguished section $v_i$, $i=1,2$,
then an affine bundle morphism $\zf:\bA_1\ra\bA_2$ is a {\it
morphism of special affine bundles} if
$$\zf_{\sv}(v_1)=v_2,$$
i.e. its linear part is a {\it morphism of special vector
bundles}. The morphism $\zf$ is a section of the natural bundle
$\bAff_M(\bA_1,\bA_2)$ whose fibers are morphisms of special
affine spaces -- the corresponding fibres. The set of morphisms of
special affine bundles $\bA_i$ will be denoted by
$\bAff(\bA_1,\bA_2)$. Similarly, bi-affine morphisms are called
{\it special bi-affine} if they are special affine with respect to
every argument separately.

In the category of special affine bundles there is a canonical
notion of duality. The {\it special affine dual} $\bA^\#$ of a
special affine bundle $\bA=(A,v_\bA)$ is an affine subbundle in
$A^\dag$ that consists of those affine functions on fibers of $A$
the linear part of which maps $v_\bA$ to $1$, i.e. those that are
special affine morphisms between $\bA$ and $\bI$:
$$\bA^\#=\bAff_M(\bA, \bI).$$ The model vector bundle for the dual
$\bA^\#$ is the vector subbundle of $\bA^\dag$ of those functions
whose linear parts vanish on $v_\bA$. We see therefore that $1_A$
is a section of the model vector bundle, since its linear part is
identically $0$. Therefore we can consider $\bA^\#$ as a special
affine bundle.

Given a special affine bundle $\bA=(A,v_\bA)$, the distinguished
section $v_\bA$ gives rise to an $\R$-action on $A$, $a_x\mapsto
a_x+tv_\bA(x)$, whose fundamental vector field is
$\chi_\bA(a_x)=\frac{\xd}{\xd t}(a_x-tv_\bA(x))$. Of course,
$\chi_\bA$ determines the section $v_\bA$.

The special affine duality is a true duality: in finite dimensions
we have a canonical identification $$(\bA^\#)^\#\simeq \bA$$ and
the corresponding special bi-affine pairing
$$\la\cdot,\cdot\ran_{sa}:\bA\ti\bA^\#\ra\bI.$$
In this identification sections of $\bA$ correspond to affine
functions $F$ on $\bA^\#$ such that $\chi_{\bA^\#}(F)=-1$.

\subsection{AV-bundles and affine phase bundles} In the standard differential
geometry many constructions are based on the algebra $C^\infty(M)$
of smooth functions on the manifold $M$. In the geometry of affine
values (AV-geometry in short) we replace $C^\infty(M)$ by the
space of sections of certain affine bundle over $M$. Let $\zz:
\Z\ra M$ be a one-dimensional affine bundle over the manifold $M$
modelled on the trivial special vector bundle $M\times \R$. Such a
bundle will be called a {\it bundle of affine values} ({\it
AV-bundle} in short). In other words, AV-bundles are
one-dimensional special affine bundles. Every special affine
bundle $\bA=(A, v_{\bA})$ gives rise to a certain bundle of affine
values $\Z=\bAP(\bA)$. The total space of the bundle $\bAP(\bA)$
is $A$ and the base is $A/\la v_\bA\ran$ (that will be denoted by
$\ul{\bA}$) with the canonical projection from the space to the
quotient. The meaning of the quotient is obvious: the class
$[a_p]$ of $a_p\in A_p$ is the orbit  $\{ a_p+tv_\bA(p)\}$ of
$\chi_\bA$. For the reason of a further application we choose the
distinguished section $v_{\bAP(\bA)}$ of the model vector bundle
characterized by $\chi_{\bAP(\bA)}=-\chi_{\bA}$.

Since $\Z$ is modelled on the trivial bundle $M\times \R$, there
is a free and transitive action of $\R$ in every fiber of $\Z$.
Therefore $\Z$ can be equivalently viewed as a principal
$\R$-bundle with the fundamental vector field induced by the
$\R$-action being $\chi_\Z$. The affine space $\Sec(\Z)$ is
modelled on the space $C^\infty(M)$ of smooth functions on the
base manifold $M$. Any section $\sigma$ of $\Z$ induces a
trivialization
    $$ I_\sigma : \Z\ni z\longmapsto  (\zz(z), z-\sigma(\zz(z)))
    \in M\times\R.$$
    The above trivialization induces the identification between
$\Sec(\Z)$ and $C^\infty(M)$. We can go surprisingly far in many
constructions replacing the ring of smooth functions with the
affine space of sections of an AV-bundle. It is possible because
many objects of standard differential geometry (like the cotangent
bundle with its canonical symplectic form) have properties that
are conserved by certain affine transformations. To build an
AV-analog of the cotangent bundle $\sT^\s M$, let us define an
equivalence relation in the set of all pairs $(x,\zs )$, where $x$
is a point in $M$ and $\sigma $ is a section of $\zz $. Two pairs
$(x,\zs )$ and $(x',\zs ')$ are equivalent if $x' = x$ and
$\xd(\zs ' - \zs)(x) = 0$.  We have identified the section $\zs '
- \zs $ of $pr_M:M\ti\R\ra M$ with a function on $M$ for the
purpose of evaluating the differential $\xd(\zs ' - \zs)(x)$. We
denote by $\sP\Z$ the set of equivalence classes. The class of
$(x,\zs )$ will be denoted by $\bd\sigma (x)$  and will be called
the {\it differential} of $\zs $ at $x$. We will write $\bd$ for
the affine exterior differential to distinguish it from the
standard  $\xd$. We define a mapping ${\sP}\zeta \colon \sP \Z
\rightarrow M$ by $\sP\zz (\bd\zs (x)) = x$. The bundle $\sP\zz$
is canonically an affine bundle modelled on the cotangent bundle
$\zp _M \colon \sT^{\textstyle *} M \rightarrow M$ with the affine
structure
$$\bd\zs _2(x)- \bd\zs _1(x) =
\xd(\zs _2 - \zs _1)(x).
$$
This affine bundle is called the {\it phase bundle} of $\zz$. A
section of $\sP\zz$ will be called an {\it affine 1-form}.

Like the cotangent bundle $\sT^\ast M$ itself, its AV-analog
$\sP\Z$ is equipped with a canonical symplectic structure. For a
chosen section $\zs$ of $\zz$ we have the isomorphisms
      \bea\nn
   &I_\zs\colon\Z\rightarrow M\times \R,  \\
      &I_{\bd\zs}\colon \sP\Z\rightarrow \sT^{\textstyle *} M,\nn
                                       \eea
   and for two sections $\zs, \zs'$ the mappings $I_{\bd\zs}$ and
   $I_{\bd\zs'}$ differ by the translation by $\xd(\zs -\zs')$, i.e.
      \bea\nn
   &I_{\bd\zs'}\circ I_{\bd\zs}^{-1}\colon \sT^{\textstyle *} M\rightarrow
\sT^{\textstyle *} M  \\
   &\colon\za_x \mapsto \za_x +\xd(\zs -\zs')(x).\nn
                                       \eea
Now we use an affine property of the canonical symplectic form
$\zw_M$ on the cotangent bundle: {\bf translations in $\sT^\s M$
by closed 1-forms are symplectomorphisms}, to conclude that the
two-form $I_{\bd\zs}^{\textstyle *} \zw_M$, where $\zw_M$ is the
canonical symplectic form on $\sT^{\textstyle *} M$, does not
depend on the choice of $\zs$ and therefore it is a canonical
symplectic form on $\sP\Z$. We will denote this form by $\zw_\bZ$.

Given a special affine bundle $\bA$, the phase bundle
$\sP(\bAP(\bA))$ will be denoted shortly $\sP\bA$.

\subsection{Brackets on affine bundles} An {\it affgebra} is an
affine space $A$ with a bi-affine operation
$$[\cdot,\cdot]: A\times A\ra \sV(A).$$ A {\it Lie affgebra } is
an affgebra with the operation satisfying the following
conditions:
  \begin{itemize}
   \item skew-symmetry: $[a_1,a_2]=-[a_2,a_1]$
   \item Jacobi identity:
   $[a_1,[a_2,a_3]]_{\sv}^2+[a_2,[a_3,a_1]]_{\sv}^2+[a_3,[a_1,a_2]]_{\sv}^2=0$.
      \end{itemize}

In the above formula $[\cdot,\cdot]_{\sv}^2$ denotes the
affine-linear part of the bracket $[\cdot,\cdot]$, i.e. we take
the linear part with respect to the second argument. One can show
that any skew-symmetric bi-affine bracket as above is completely
determined by its affine-linear part. The detailed description of
the concept of Lie affgebra can be found in \cite{GGU1}.

A {\it Lie affgebroid} is an affine bundle $A$ over $M$ with a Lie
affgebra bracket on the space of sections
$$[\cdot,\cdot]: \Sec(A)\times \Sec(A)\ra \Sec(\sV(A))$$
together with a morphism of affine bundles $\rho:A\ra\sT M$
inducing a map from $\Sec(A)$ into vector fields on $M$ such that
$$[a,fv]_{\sv}^2=f[a,v]_{\sv}^2+\rho(a)(f)v$$
for any smooth function $f$ on $M$. Note that the same concept has
been introduced by E.~Mart\'\i nez, T.~Mestdag and W.~Sarlet under
the name of {\it affine Lie algebroid} \cite{MMS1}.

A {\it special Lie affgebroid} is a special affine bundle
$\bA=(A,v_\bA)$ equipped with a Lie affgebroid bracket
$[\cdot,\cdot]$ such that $v_\bA$ is central, i.e.
$[a,v_A]_{\sv}^2=0$ for all $a\in\Sec(A)$. It is easy to see that
the special Lie affgebroid bracket on $\bA$ induces a Lie
affgebroid bracket on $\ul{\bA}$.

\begin{example} Given a fibration $\zx\colon M\ra\R$, take the affine
subbundle $A\subset\sT M$ characterized by  $\zx_*(X)=\pa_t$ for
$X\in\Sec(A)$. Then the standard bracket of vector fields in $A$,
and $\zr\colon A\ra\sT M$ being just the inclusion, define on $A$
a structure of a Lie affgebroid. This is the basic example of the
concept of {\it affine Lie algebroid} developed in
\cite{MMS1,MMS}. The special vector bundle $\cT M=\sT M\ti\R$,
whose sections represent first-order differential operators
$(X+f)$ on $M$ and the distinguished section is represented by the
constant function $f=1$,  is canonically a {\it special Lie
algebroid}, i.e. it is a Lie algebroid with the bracket
$[X+f,Y+g]=[X,Y]_{vf}+X(g)-Y(f)$ and the distinguished section $1$
is a central element with respect to this bracket. Here clearly
$[\cdot,\cdot]_{vf}$ is the bracket of vector fields. Now, the
affine subbundle $A\ti\R$ in $\cT M$ whose sections $X+f$ satisfy
$\zx_*(X)=\pa_t$ is a special Lie affgebroid.
\end{example}

Let us remind that any affine bundle $A$ is canonically embedded
in the vector bundle $\widehat A$ being its vector hull. It turns
out that there is a one-to-one correspondence between Lie
affgebroid structures on $A$ and  Lie algebroid structures on the
vector hull (see \cite[Theorem 11]{GGU1}). For special Lie
affgebroids the correspondence is analogous.

\begin{theorem} For a special affine bundle $\bA$ the following are
equivalent:
\begin{description}
\item{(a)} The bracket $[\cdot,\cdot]$ is a special Lie affgebroid
bracket on $\bA$; \item{(b)} The bracket $[\cdot,\cdot]$ is the
restriction of a special Lie algebroid bracket $[\cdot,\cdot]^\we$
on the special vector bundle $\wh{\bA}$.
\end{description}
\end{theorem}


\begin{definition} Given an AV-bundle $\Z$ over $M$, an {\it aff-Poisson
bracket} on $\Z$ is a Lie affgebra bracket
$$\{\cdot,\cdot\}\colon\Sec(\Z)\ti\Sec(\Z)\ra C^\infty(M)$$
such that
$$X_\zs=\{\zs,\cdot\}^2_\sv\colon C^\infty(M)\ra
C^\infty(M)$$ is a a vector field on $M$ (called {\it the
Hamiltonian vector field of $\zs$}) for every $\zs\in\Sec(\Z)$.
\end{definition}

For any $\bAP$-bundle $\bZ=(Z,v_\Z)$ the tangent bundle $\sT \bZ$
is equipped with the tangent $\R$-action. Dividing $\sT\bZ$ by the
action we obtain the Atiyah algebroid of the principal $\R$-bundle
$\Z$ which we denote by $\widetilde\sT \bZ$. It is a special Lie
algebroid whose sections are interpreted as invariant vector field
on the principal $\R$-bundle $\Z$. The distinguished section of
$\widetilde\sT \bZ$ is represented by the fundamental vector field
$\chi_\bZ$. The AV-bundle $\bAP(\widetilde\sT \bZ)$ is a bundle
over $\sT M$. The special affine dual for the special affine
structure on $\widetilde\sT \bZ$ is $(\widetilde\sT \bZ)^\#=\sP
\bZ\times \bI$. The special affine evaluation between
$\sP\bZ\times \bI$ and $\widetilde\sT \bZ$  comes from the
interpretation of sections of $\sP\Z$ as $\R$-invariant 1-forms
$\za$ on $\Z$ such that $\la\za,\chi_\bZ\ran=1$ (i.e. principal
connections on the $\R$-principal bundle $\Z$) that gives an
affine-linear pairing
$\la\cdot,\cdot\ran_\dag:\sP\Z\ti_M\wt{\sT}\Z\ra\R$ and the
identification $(\sP\Z)^\dag=\wt{\sT}\Z$ (cf. \cite{GGU1,GGU2}).
In this way sections of $(\sP\Z)^\dag=\wt{\sT}\Z$ represent {\it
affine derivations} $D:\Sec(\Z)\ra C^\infty(M)$ (i.e. such affine
maps that $D_\sv:C^\infty(M)\ra C^\infty(M)$ is a derivation, thus
a vector field on $M$). It is now obvious that affine
biderivations
$$\{\cdot,\cdot\}\colon\Sec(\Z)\ti\Sec(\Z)\ra C^\infty(M)$$ are
sections of the bundle $\wt{\sT}\Z\ot_M\wt{\sT}\Z$. In this
picture, skew-symmetric affine biderivations are sections of
$\wedge^2\wt{\sT}\Z$.
\subsection{Local affine coordinates and canonical identifications}
For a special affine bundle $\bA=(A,v_\bA)$ of rank $m$, we use a
section $e_0$ of the affine bundle $A$ to identify $A$ with its
model special vector bundle $\sV(\bA)$. The distinguished section
$v_\bA$ we can then extend to a basis $e_1,\dots,e_m$ of local
sections of $\sV(\bA)$ such that $e_m=v_\bA$. In this way we can
get a basis $e_0,\dots,e_m$ of local sections of $\wh{\bA}$ and
the dual basis $e^\*_0,\dots,e^\*_m$ of local sections in
$\wh{\bA}^\*=\wh{\bA^\#}$ such that $e^\*_0$ represents the
distinguished section of $\sV(\bA^\#)$. If we choose local
coordinates $(x^a)$ on the base manifold, these bases give rise to
local coordinates $(x^a,y^0,\dots,y^m)$ and
$(x^a,\zx_0,\dots,\zx_m)$ in the special vector bundles $\wh{\bA}$
and $\wh{\bA}^\*$, respectively, defined by
\bea\label{coo}y^i&=&\zi_{e^\*_i}\quad\text{for}\quad i=0,\dots,m-1,
\qquad y^m=-\zi_{e^\*_m},\\
\zx_i&=&\zi_{e_i}\quad\text{for}\quad i=1,\dots,m, \qquad
\zx_0=-\zi_{e_0}.\label{coo1}
\eea
With these coordinates the affine subbundles $\bA$ and $\bA^\#$ in
$\wh{\bA}$ and $\wh{\bA^\#}$ are characterized by the equations
$y^0=-1$ and $\zx_m=-1$, respectively. Moreover
$(x^a,y^1,\dots,y^m)$ and $(x^a,\zx_0,\dots,\zx_{m-1})$ are local
coordinates in $\bA$ and $\bA^\#$, respectively, in which the
special affine pairing $\la\cdot,\cdot\ran_{sa}:\bA\ti\bA^\#\ra\R$
reads
\be\label{par}\la(x,y^1,\dots,y^m),(x,\zx_0,\dots,\zx_{m-1})\ran_{sa}=
\sum_1^{m-1}y^i\zx_i-y^m-\zx_0.\ee Note that the base manifold
$\ul{\bA}$ of the AV-bundle $\zz:\bAP(\bA)\ra \ul{\bA}$ is an
affine bundle $\ul{\zh}:\ul{\bA}\ra M$ with induced coordinates
$(x^a,y^i)$, $i=1,\dots,m-1$, so
$\zz(x,y^1,\dots,y^m)=(x,y^1,\dots,y^{m-1})$ and
$\chi_{\bAP(\bA)}=-\pa_{y^m}$, i.e.
$$\bA\ni(x,y^1,\dots,y^m)\mapsto((x,y^1,\dots,y^{m-1}),y^m)\in\ul{\bA}\ti\bI$$
represents a local isomorphism of $\bAP(\bA)$ and
$\ul{\bA}\ti\bI$. Similar observations hold for $\bA^\#$ and
coordinates $(x,\zx)$. Coordinates on special affine bundles used
in the sequel will be always of this type.

In the space of sections of $\bAP(\bA)$ one can distinguish {\it
affine sections}, i.e. such sections
$$\zs:\ul{\bA}=A/\la v_\bA\ran\ra \bA$$
which are affine morphisms. The space of affine sections will be
denoted by $\Aff\Sec(\bAP(\bA))$. We say that an operation
$\Sec(\bAP(\bA))\ti\Sec(\bAP(\bA))\ra C^\infty(\ul{\bA})$ is {\it
affine-closed} if the product of any two affine sections is an
affine function.

In the linear case there is a correspondence between Lie algebroid
brackets on a vector bundle and linear Poisson structures on the
dual bundle. In the affine setting we have an analog of this
correspondence. Let $X$ be a section of $\sV(\bA)$. Since
$\sV(\bA)$ is a vector subbundle in $\wh \bA$, the section X
corresponds to a linear function $\iota^\dag_X$ on $\bA^\dag=(\wh
\bA)^\*$. The function $\iota^\dag_X$ is invariant with respect to
the vertical lift of the distinguished section $1_A$ of
$\bA^\dag$, so its restriction to $\bA^\#$ is constant on fibres
of the projection $\bA^\#\ra\bA^\#/\la 1_A\ran$ and defines an
affine function $\iota^\#_X$ on the base
$\ul{\bA^\#}=\bA^\#\slash\langle 1_A\rangle$ of $\bAP(\bA^\#)$.
Hence we have a canonical identification between
\begin{description}
\item{(a)}  sections $X$ of $\sV(\bA)$, \item{(b)} linear
functions $\iota_X^\dag$ on $\bA^\dag$ which are invariant with
respect to the vertical lift of $1_A$, \item{(c)} affine functions
$\iota_X^\#$ on $\ul{\bA^\#}$.
\end{description}
In local coordinates and bases, if $X=\sum_1^mf_i(x)e_i$, then
$\iota_X^\dag=\sum_1^{m-1}f_i(x)\zx_i-f_m(x)\zx_m$ and
$\iota_X^\#=f_m(x)+\sum_1^{m-1}f_i(x)\zx_i$.

On the other hand, in the theory of vector bundles there is an
obvious identification between sections of the bundle $E$ and
functions on $E^\ast$ which are linear along fibres. If $\zf$ is a
section of $E$, then the corresponding function $\zi_\zf$ is
defined by the canonical pairing
$$\zi_\zf(X)=\langle \zf,X\rangle.$$
In the theory of special affine bundles we have an analog of the
above identification:
$$\Sec(\bA)\simeq\Aff\Sec(\bAP(\bA^\#)),\quad
\bF_\zs\leftrightarrow\zs,$$ where  $\bF_\zs$ is the unique
(affine) function on $\bA^\#$ such that
$\chi_{\bAP(\bA^\#)}(\bF_\zs)=1$ ($\bF_\zs$ can be therefore
interpreted as a section of $\bA$) and $\bF_\zs\circ\zs=0$. In
local affine coordinates, we associate with the section
$$\zx_0=\zs(x^a,\zx_1,\dots,\zx_{m-1})=\sum_1^{m-1}\zs_i(x)\zx_i-\zs_m(x)$$
of $\bAP(\bA^\#)$ the function
$$\bF_\zs=\zs(x^a,\zx_1,\dots,\zx_{m-1})-\zx_0=
\sum_1^{m-1}\zs_i(x)\zx_i-\zs_m(x)-\zx_0$$ on $\bA^\#$ which
represents the section $M\ni x\mapsto\sum_1^{m}\zs_i(x)e_i\in\bA$.

\medskip
Using the above identification we can formulate the following
theorem (for the proof see \cite{GGU2}):

\begin{theorem}
There is a canonical one-to-one correspondence between special Lie
affgebroid brackets $[\cdot,\cdot]_\bA$ on $\bA$ and affine
aff-Poisson brackets $\{\cdot,\cdot\}_{\bA^\#}$ on $\bAP(\bA^\#)$,
uniquely defined by:
$$\{\zs,\zs'\}_{\bA^\#}=\zi^\#_{[\bF_\zs,\bF_{\zs'}]_\bA}.$$
\end{theorem}

\begin{example} Let $\bZ=(Z,v_\bZ)$ be an $\bAP$-bundle.
The AV-bundle $\bAP((\wt{\sT}\Z)^\#)$ is the trivial bundle over
the affine phase bundle $\sP\Z$ and the aff-Poisson bracket,
associated with the canonical special Lie algebroid bracket on
$\wt{\sT}\Z$, is the standard Poisson bracket on $\sP\Z$
associated with the canonical symplectic form $\zw_\Z$ on $\sP\Z$.
\end{example}

\section{Double bundles}
Let $M$ be a smooth manifold and let $(x^a), \ a=1,\dots,n$, be a
coordinate system in $M$. We denote by $\zt_M \colon \sT M
\rightarrow M$ the tangent vector bundle and by $\zp_M \colon
\sT^\* M\rightarrow M$ the cotangent vector bundle. We have the
induced (adapted) coordinate systems $(x^a, {\dot x}^b)$ in $\sT
M$ and $(x^a, p_b)$ in $\sT^\* M$.
        Let $\zt\colon E \rightarrow M$ be a vector bundle and let $\zp
\colon E^\* \rightarrow M$ be the dual bundle.
  Let $(e_1,\dots,e_m)$  be a basis of local sections of $\zt\colon
E\rightarrow M$ and let $(e^{1}_*,\dots, e^{m}_*)$ be the dual
basis of local sections of $\zp\colon E^\*\rightarrow M$. We have
the induced coordinate systems:
    \beas
    (x^a, y^i),\quad & y^i=\zi(e^{i}_*), \quad \text{in} \ E,\\
    (x^a, \zx_i), \quad &\zx_i = \zi(e_i),\quad \text{in} \ E^\* ,
    \eeas
    where the linear functions  $\zi(e)$ are given by the canonical pairing
    $\zi(e)(v_x)=\la e(x),v_x\ran$. Thus we have local coordinates
    \beas
    (x^a, y^i,{\dot x}^b, {\dot y}^j ) &  \quad \text{in} \ \sT E,\\
    (x^a, \zx_i, {\dot x}^b, {\dot \zx}_j) & \quad \text{in} \ \sT E^\* ,\\
    (x^a, y^i, p_b, \zp_j) & \quad \text{in}\ \sT^\*E,\\
    (x^a, \zx_i, p_b, \zf^j) & \quad \text{in}\ \sT^\* E^\* .
    \eeas
It is well known (cf. \cite{KU}) that the cotangent bundles
$\sT^\*E$ and $\sT^\*E^\*$ are examples of double vector bundles:
$$\xymatrix{
\sT^\ast E^\ast\ar[rr]^{\sT^\ast\zp} \ar[d]_{\zp_{E^\ast}} && E\ar[d]^{\zt} \\
E^\ast\ar[rr]^{\zp} && M } \qquad {,}\qquad \xymatrix{
\sT^\ast E\ar[rr]^{\sT^\ast\zt} \ar[d]_{\zt_{E^\ast}} && E^\ast\ar[d]^{\zp} \\
E\ar[rr]^{\zt} && M }.
$$
Note that the concept of a double vector bundle is due to
J.~Pradines \cite{Pr1,Pr2}, see also \cite{Ma,KU}. In particular,
all arrows correspond to vector bundle structures and all pairs of
vertical and horizontal arrows are vector bundle morphisms. The
above double vector bundles are canonically isomorphic with the
isomorphism
    $$\cR_\zt \colon \sT^\*E \longrightarrow \sT^\* E^\*
                    $$
  being simultaneously an anti-symplectomorphism  (cf. \cite{Du, KU, GU2}):
\be\label{iso}\xymatrix{
 & \sT^\ast E^\ast  \ar[dr]^{\sT^\ast\pi}
 \ar[ddl]_{\pi_{E^\ast}}
 & & & \sT^\ast E\ar[dr]^{\pi_E}\ar[ddl]_/-20pt/{\sT^\ast\tau}\ar[lll]_{\cR_r}
 & \\
 & & E\ar[ddl]_/-20pt/{\tau}
 & & & E \ar[ddl]_{\tau}\ar[lll]_/20pt/{id}\\
 E^\ast\ar[dr]^{\pi}
 & & & E^\ast\ar[dr]^{\pi}\ar[lll]_/20pt/{id} & &  \\
 & M& & & M\ar[lll]_{id} &
}\ee In local coordinates, $\cR_\zt$ is given by
    $$\cR_\zt(x^a, y^i, p_b, \zp_j) = (x^a, \zp_i, -p_b,y^j).
                              $$

Exactly like the cotangent bundle $\sT^\*E$ is a double vector
bundle, the affine phase bundle $\sP\bA$ is canonically a {\it
double affine bundle}:
$$\xymatrix{
\sP\bA\ar[rr]^{\sP^\#\zz} \ar[d]_{\sP\zz} && \ul{\bA^\#}\ar[d]^{\ul{\zh^\#}} \\
\ul{\bA}\ar[rr]^{\ul{\zh}} && M }.$$

\medskip\noindent
{\bf Definition.} \begin{enumerate} \item{} A {\it trivial double
affine bundle} over a base manifold $M$ is a commutative diagram
of four affine bundle projections
\be\label{dab}\xymatrix{
A\ar[r]^{\zX_2} \ar[d]_{\zX_1} & A_2\ar[d]^{\zh_2} \\
A_1\ar[r]^{\zh_1} & M },\ee where $A=M\ti K_1\ti K_2\ti K$ is a
trivial affine bundle with the fibers being direct products of
(finite-dimensional) affine spaces $K_1,K_2,K$, with the obvious
projections:
$$\zX_i:M\ti K_1\ti K_2\ti K\ra A_i=M\ti K_i,\quad
\zh_i:A_i=M\ti K_i\ra M,\quad i=1,2.$$ In particular, the pairs
$(\zX_i,\zh_i)$, $i=1,2$, are morphism of affine bundles. \item A
{\it morphism of trivial double affine bundles over the identity
on the base} is a commutative diagram of morphisms of trivial
affine bundles as above
\be\label{dam}\xymatrix{
 & A \ar[dr]^{\zX_2}
 \ar[ddl]_{\zX_1}
 & & & A'\ar[dr]^{\zX'_2}\ar[ddl]_/-20pt/{\zX'_1}\ar[lll]_{\zF}
 & \\
 & & A_2\ar[ddl]_/-20pt/{{\zh_2}}
 & & & A'_2 \ar[ddl]_{{\zh'_2}}\ar[lll]_/20pt/{\zF_2}\\
 A_1\ar[dr]^{\zh_1}
 & & & A_1'\ar[dr]^{\zh'_1}\ar[lll]_/20pt/{\zF_1} & &  \\
 & M& & & M\ar[lll]_{id} &
}\ee where clearly $A=M\ti K_1\ti K_2\ti K$,  $\ A'=M\ti K'_1\ti
K'_2\ti K'$, etc., and $(\zX_i,\zh_i)$, $(\zX'_i,\zh'_i)$,
$(\zF,\zF_i)$, and $(\zF_i, id)$, $i=1,2$, are morphisms of affine
bundles. Note that we have not considered $A$ as an affine bundle
over $M$, since this structure is accidental in the trivial case
and it is not respected by isomorphisms in the above sense.

\item A {\it double affine bundle} modelled on the product $K_1\ti
K_2\ti K$ of affine spaces is a commutative diagram (\ref{dab}) of
affine bundles, this time not necessarily trivial, which is
locally diffeomorphic with trivial double affine bundles
associated with $K_1\ti K_2\ti K$, i.e. there is an open covering
$\{ U_\za\}$ of $M$ such that
$$\xymatrix{
A^\za\ar[rr]^{{\zX_2}_{\mid A^\za}} \ar[d]_{{\zX_1}_{\mid A^\za}}
&&
A^\za_2\ar[d]^{{\zh_2}_{\mid A^\za_2}} \\
A^\za_1\ar[rr]^{{\zh_1}_{\mid A^\za_1}} && M },$$ where
$A_i^\za=\zh_i^{-1}(U_\za)$, $i=1,2$, and
$A^\za=(\zX_1\circ\zh_1)^{-1}(U_\za)=(\zX_2\circ\zh_2)^{-1}(U_\za)$,
is (diffeomorphically) equivalent to the trivial double affine
bundle
$$\xymatrix{
U_\za\ti K_1\ti K_2\ti K\ar[rr] \ar[d] &&
U_\za\ti K_2\ar[d] \\
U_\za\ti K_1\ar[rr] && M }$$ and such that this equivalence
induces an automorphism of the corresponding trivial double affine
bundle over the identity on each intersection $U_\za\bigcap
U_\zb$.

\item A {\it morphism of double affine bundles} is a commutative
diagram (\ref{dam}) of morphisms of corresponding affine bundles
which in local trivializations induces morphisms od trivial double
affine bundles
\end{enumerate}

\medskip\noindent
To see how the trivial double affine bundles are glued up inside a
(non-trivial) double affine bundle, let us take affine coordinates
$(x,y,z,c)$ in $A=\R^k\ti K_1\ti K_2\ti K$. Then every morphism
$\zF:A\ra A$ of the double affine bundle into itself (over the
identity on the base) is of the form
\beas \zF(x,y^j,z^a,c^u)&=&(x,\,\za_0^j(x)+\sum_i\za^j_i(x)y^i,\,
\zb_0^a(x)+\sum_b\zb_b^a(x)z^b,\\
&&\zg_{00}^u(x)+\sum_i\zg_{i0}^u(x)y^i+\sum_b\zg_{0b}^u(x)z^b+
\sum_{i,b}\zg_{ib}^u(x)y^iz^b+\sum_w\zs^u_w(x)c^w).
\eeas

\medskip\noindent
Let $\bA$ be a special affine bundle. The adapted affine
coordinates on $\sP\bA$ are $(x^a,y^i,p_b,\zx_j)$,
$i,j=1,\dots,m-1$, and the projections read
$$\sP^\#\zz(x,y,p,\zx)=(x,\zx),\qquad
\ul{\zh^\#}(x,y,p,\zx)=(x,y).$$ The canonical symplectic form
reads $\zw_{\sP\bA}=\xd p_a\we\xd x^a+\xd\zx_i\we\xd y^i$ and the
affine de Rham differential $\bd:\Sec(\bAP(\bA))\ra\Sec(\sP\bA)$
associates with a section $y^m=\zs(x,y^i)$, $i=1,\dots,m-1$ the
section $\bd\zs$:
$$p_a\circ\bd\zs=\frac{\pa \zs}{\pa x^a},\quad
\zx_i\circ\bd\zs=\frac{\pa \zs}{\pa y^i}.$$ The special vector
bundle $\wt{\sT}(\bAP(\bA^\#))=(\sP\bA^\#)^\dag$, which be denoted
shortly $\wt{\sT}\bA^\#$, carries the coordinates
$(x^a,y^i,\dot{x}^b,\dot{y}^j,s)$, $i,j=1,\dots,m-1$, with the
affine-linear pairing
$\la\cdot,\cdot\ran_\dag:\sP\bA^\#\ti_{\ul{\bA^\#}}\wt{\sT}\bA^\#\ra\R$
given by
$$\la(x,y,\zx),(x,y,\dot{y},s)\ran_\dag=\zx_i\dot{y}^i-s.$$
We have the obvious projection $\wt{\sT}\bA^\#\ra\sT\ul{\bA^\#}$
that reads $(x,y,\dot{x},\dot{y},s)\mapsto(x,y,\dot{x},\dot{y})$.
Similarly like in the linear case, the double affine bundles
$\sP\bA$ and $\sP\bA^\#$ are canonically isomorphic (see
\cite{Ur1})
\be\label{iso1}\xymatrix{
 & \sP\bA^\#  \ar[dr]^{\sP^\#\zh^\#}
 \ar[ddl]_{\sP\zh^\#}
 & & & \sP\bA\ar[dr]^{\sP\zh}\ar[ddl]_/-20pt/{\sP^\#\zh}\ar[lll]_{\cR_\zh}
 & \\
 & & \ul{\bA}\ar[ddl]_/-20pt/{\ul{\zh}}
 & & & \ul{\bA} \ar[ddl]_{\ul{\zh}}\ar[lll]_/20pt/{id}\\
 \ul{\bA^\#}\ar[dr]^{\ul{\zh^\#}}
 & & & \ul{\bA^\#}\ar[dr]^{\ul{\zh^\#}}\ar[lll]_/20pt/{id} & &  \\
 & M& & & M\ar[lll]_{id} &
}\ee In local coordinates, $\cR_\zh$ is given by
    $$ \cR_\zh(x^a,y^i,p_b,\zx_j) = (x^a, \zx_i, -p_b,y^j).
                              $$

\section{Lagrangian and Hamiltonian formalisms \\ for general
algebroids}

\subsection{Lie algebroids as double vector bundle morphisms\label{S1}}
For Lie algebroids we refer to the survey article \cite{Ma1}. It
is well known that Lie algebroid structures on a vector bundle $E$
correspond to linear Poisson tensors on $E^\*$. A 2-contravariant
tensor $\zL$ on $E^\*$ is called {\it linear} if the corresponding
mapping $\widetilde{\zL} \colon \sT^\* E^\* \rightarrow \sT E^\*$
induced by contraction is a morphism of double vector bundles.
This is the same as to say that the corresponding bracket of
functions is closed on (fiber-wise) linear functions. The
commutative diagram
$$\xymatrix{
\sT^\ast E^\ast\ar[r]^{\widetilde\Lambda}  & \sT E^\ast \\
\sT^\ast E\ar[u]_{\cR_\tau}\ar[ur]^{\ze} & },
$$
composed with (\ref{iso}), describes a one-to-one correspondence
between linear 2-contravariant tensors $\zL$ on $E^\*$ and
homomorphisms of double vector bundles (cf. \cite{KU, GU2})
covering the identity on $E^\*$:

\be\xymatrix{
 & \sT^\ast E \ar[rrr]^{\varepsilon} \ar[dr]^{\pi_E}
 \ar[ddl]_{\sT^\ast\tau}
 & & & \sT E^\ast\ar[dr]^{\sT\pi}\ar[ddl]_/-20pt/{\tau_{E^\ast}}
 & \\
 & & E\ar[rrr]^/-20pt/{\varepsilon_l}\ar[ddl]_/-20pt/{\tau}
 & & & \sT M \ar[ddl]_{\tau_M}\\
 E^\ast\ar[rrr]^/-20pt/{id}\ar[dr]^{\pi}
 & & & E^\ast\ar[dr]^{\pi} & &  \\
 & M\ar[rrr]^{id}& & & M &
}\label{F1.3}
\ee
In local coordinates, every  $\ze$ as in above is of the form
\be\label{F1.4}
\ze(x^a,y^i,p_b,\zp_j) = (x^a, \zp_i, \sum_k\zr^b_k(x)y^k,
\sum_{i,k}c^k_{ij}(x) y^i\zp_k + \sum_a\zs^a_j(x) p_a)
\ee
and it corresponds to the linear tensor
$$ \zL_\ze =\sum_{i,j,k}c^k_{ij}(x)\zx_k \partial _{\zx_i}\otimes \partial _{\zx_j} +
\sum_{i,b}\zr^b_i(x) \partial _{\zx_i} \otimes \partial _{x^b} -
\sum_{a,j}\zs^a_j(x)\partial _{x^a} \otimes \partial _{\zx_j}.
$$
In \cite{GU2} by  {\it algebroids} we meant the morphisms
(\ref{F1.3}) of double vector bundles covering the identity on
$E^\*$, while {\it Lie algebroids} were those algebroids for which
the tensor $\zL_\ze$ is a Poisson tensor. The relation to the
canonical definition of Lie algebroid is given by the following
theorem (cf. { \cite{GU3, GU2}}).

\begin{theo}
An algebroid structure $(E,\ze)$ can be equivalently defined as a
bilinear bracket $[\cdot ,\cdot]_\ze $ on sections of $\zt\colon
E\rightarrow M$, together with vector bundle morphisms $\ze_l,\,
\ze_r \colon E\rightarrow \sT M$ (left and right anchors), such
that
$$ [fX,gY]_\ze = f(\ze_l\circ X)(g)Y -g(\ze_r \circ Y)(f) X
+fg [X,Y]_\ze
$$
         for $f,g \in \cC^\infty (M)$, $X,Y\in \otimes ^1(\zt)$.
The bracket and anchors are related to the 2-contravariant tensor
$\zL_\ze$ by the formulae
\beas
        \zi([X,Y]_\ze)&= \{\zi(X), \zi(Y)\}_{\zL_\ze},  \\
        \zp^\*(\ze_l\circ X(f))       &= \{\zi(X), \zp^\*f\}_{\zL_\ze}, \\
        \zp^\*(\ze_r\circ X(f))       &= \{\zp^\* f, \zi(X)\}_{\zL_\ze}.
                                                   \eeas
        The algebroid $(E,\ze)$ is a Lie algebroid if and only if the tensor
$\zL_\ze$ is a Poisson tensor.
\end{theo}

The canonical example of a mapping $\ze$ in the case of $E=\sT M$
is given by $\ze = \ze_M = \za^{-1}_M$ -- the inverse to the
Tulczyjew isomorphism $\za_M$ that can be defined as the dual to
the isomorphism of double vector bundles
$$\xymatrix{
 & \sT \sT M \ar[rrr]^{\zk_M} \ar[dr]^{\zt_{\ssT M}}
 \ar[ddl]_{\sT\zt_M}
 & & & \sT \sT M\ar[dr]^{\sT\zt_M}\ar[ddl]_/-20pt/{\zt_{\ssT M}}
 & \\
 & & \sT M\ar[rrr]^/-20pt/{id}\ar[ddl]_/-20pt/{\zt_M}
 & & & \sT M \ar[ddl]_{\zt_M}\\
 \sT M\ar[rrr]^/-20pt/{id}\ar[dr]^{\zt_M}
 & & & \sT M\ar[dr]^{\zt_M} & &  \\
 & M\ar[rrr]^{id}& & & M &
}$$
In general, the algebroid structure map $\ze$ is not an isomorphism and,
consequently, its dual $\zk^{-1} = \ze^{\s_r}$  with respect to
the right projection is a relation and not a mapping.

\subsection{Lagrangian and Hamiltonian formalisms}
The double vector bundle morphism (\ref{F1.3}) can be extended to
the following algebroid analogue of the so called Tulczyjew triple
\be\xymatrix@C-5pt{
 &\sT^\ast E^\ast \ar[rrr]^{\widetilde{\Lambda}}
\ar[ddl]_{\pi_{E^\ast}} \ar[dr]^{\sT^\ast\pi}
 &  &  & \sT E^\ast \ar[ddl]_/-25pt/{\zt_{E^\ast}} \ar[dr]^{\sT\pi}
 &  &  & \sT^\ast E \ar[ddl]_/-25pt/{\sT^\ast\zt} \ar[dr]^{\pi_E}
\ar[lll]_{\varepsilon}
 & \\
 & & E \ar[rrr]^/-20pt/{\wt{\zL}_r}\ar[ddl]_/-20pt/{\zt}
 & & & \sT M\ar[ddl]_/-20pt/{\zt_M}
 & & & E\ar[lll]_/+20pt/{\varepsilon_l}\ar[ddl]_{\zt}
 \\
 E^\ast\ar[rrr]^/-20pt/{id} \ar[dr]^{\pi}
 & & & E^\ast\ar[dr]^{\pi}
 & & & E^\ast\ar[dr]^{\pi}\ar[lll]_/-20pt/{id}
 & & \\
 & M\ar[rrr]^{id}
 & & & M & & & M\ar[lll]_{id} &
}\label{F1.3b}\ee The left-hand side is Hamiltonian, the
right-hand side is Lagrangian, and the dynamics lives in the
middle.

Any Lagrangian function $L:E\ra\R$ defines a Lagrangian
submanifold $N=(\xd L)(E)$ in $\sT^\* E$, being the image of the
de Rham differential $\xd L$, i.e. the image of the section $\xd
L:E\ra\sT^\* E$. The further image $D=\ze(N)$ can be understood as
an implicit differential equation on $E^\*$, solutions of which
are `phase trajectories' of the system. The Lagrangian defines
also smooth maps: $L_{eg}:E\ra E^\*$ and $\wt{L}_{eg}:E\ra \sT
E^\*$ by
$$L_{eg}=\zt_{E^\*}\circ\ze\circ\xd L=\sT^\*\zt\circ\xd L$$ and
$\wt{L}_{eg}=\ze\circ\xd L$. The map $L_{eg}$  is {\it de facto}
the vertical derivative of $L$ and  is the analogue of the {\it
Legendre mapping}. We use this terminology to distinguish the
Legendre mapping  associated with $L$ from the Legendre
transformation which we understand as the passage from a
Lagrangian to a Hamiltonian generating object as explained in
\cite{TU}). The introduced ingredients produce an implicit
differential equation, this time for curves $\zg:I\ra E$.

This equation, which will be denoted by $(E_L)$,  is represented
by a subset $E_L$ of $\sT E$ being  the inverse image
$$E_L=\sT(\wt{L}_{eg})^{-1}(\sT^2E^\*)$$ of the subbundle
$\sT^2E^\*$ of holonomic vectors in $\sT\sT E^\*$ (i.e. such
vectors $w$ that $\zt_{\sT E^\*}(w) = \sT\zt_{E^\*} (w)$) with
respect to the derivative $\sT(\wt{L}_{eg}):\sT E\ra\sT\sT E^\*$
of $\wt{L}_{eg}:E\ra\sT E^\*$. The solutions are such paths
$\zg:I\ra E$ that the tangent prolongation $\sT(L_{eg}\circ\zg)$
of $L_{eg}\circ\zg$ is exactly $\wt{L}_{eg}\circ\zg$.

In local coordinates, $D$ has the parametrization by $(x^a,y^k)$
via $\wt{L}_{eg}$ in the form (cf. (\ref{F1.4}))
\be\wt{L}_{eg}(x^a,y^i)= (x^a,\frac{\partial L}{\partial
y^i}(x,y), \sum_k\zr^b_k(x)y^k, \sum_{i,k}c^k_{ij}(x)
y^i\frac{\partial L}{\partial y^k}(x,y) +\sum_a
\zs^a_j(x)\frac{\partial L}{\partial x^a}(x,y)) \label{F1.4a}\ee
and the equation $(E_L)$, for $\zg(t)=(x^a(t),y^i(t))$, reads \be
(E_L):\qquad\qquad\frac{\xd x^a}{\xd t}=\sum_k\zr^a_k(x)y^k,
\qquad \frac{\xd}{\xd t}\left(\frac{\partial L}{\partial
y^j}\right)= \sum_{i,k}c^k_{ij}(x) y^i\frac{\partial L}{\partial
y^k}(x,y) + \sum_a\zs^a_j(x)\frac{\partial L}{\partial
x^a}(x,y),\label{EL2}\ee in the full agreement with \cite{LMM,
Mar1, Mar2, We}, if only one takes into account that, for Lie
algebroids, $\zs^a_j=\zr^a_j$. As one can see from (\ref{EL2}),
the solutions are automatically admissible curves in $E$, i.e. the
velocity  $\frac{\xd }{\xd t}(\zt\circ\zg)(t)$ is $\ze_l(\zg(t))$.

\vskip10pt\noindent Note that the tensor $\zL_\ze$ gives rise also
to kind of a Hamiltonian formalism (cf. \cite{OPB}). In \cite{GU2}
and \cite{OPB} one refers to a 2-contravariant tensor as to a {\it
Leibniz structure} (this notion is completely different from {\it
Leibniz algebra} in the sense of J.-L.~Loday -- a
non-skew-symmetric analog of a Lie algebra). In the presence of
$\zL_\ze$, by the {\it hamiltonian vector field} associated with a
function $H$ on $E^\*$ we understand the contraction $\ix_{\xd
H}\zL_\ze$. Thus the question of the Hamiltonian description of
the dynamics $D$ is the question if $D$ is the image of a
Hamiltonian vector field. (Of course, one can also try to extend
such a Hamiltonian formalism to more general generating objects
like Morse families.) Every such a function $H$ we call a {\it
Hamiltonian  associated with the Lagrangian} $L$. However, it
should be stressed that, since $\ze$ and $\zL_\ze$ can be
degenerate, we have much more freedom in choosing generating
objects (Lagrangian and Hamiltonian) than in the symplectic case.
For instance, the Hamiltonian is defined not up to a constant but
up to a Casimir function of the tensor $\zL_\ze$ and for the
choice of the Lagrangian we have a similar freedom. However, in
the case of a hyperregular Lagrangian, i.e. when $L_{eg}$ is a
diffeomorphism, we recover the standard correspondence between
Lagrangians and Hamiltonians. We have namely the following (see
\cite[Corollary 1]{GGU3}).

\begin{theo} \label{ham}
If the Lagrangian $L$ is hyperregular, then the function
$$H(e^\*_x)=\la L_{eg}^{-1}(e^\*_x),e^\*_x\ran-L\circ
L_{eg}^{-1}(e^\*_x)$$ is a Hamiltonian associated with $L$. This
Hamiltonian has the property that the Lagrange submanifold $N=\xd
L(M)$ in $\sT^\* E$ corresponds under the canonical isomorphism
$\cR_\zt$ to the Lagrange submanifold $\xd H(M)$ in $\sT^\* E^\*$.
\end{theo}

\section{Lagrangian and Hamiltonian formalisms \\ for special
affgebroids}

\subsection{Special affgebroids as morphisms of double affine
bundles}

Let $\zh:\bA=(A, v_{A})\ra M$ be a special affine bundle. With an
analogy to the linear case, by a general {\it special affgebroid
on $\bA$} we mean a morphism of double affine bundles covering the
identity on $\ul{\bA^\#}$:
\be\xymatrix{
 & \sP\bA\ar[rrr]^{\cE} \ar[dr]^{\sP\zh}
 \ar[ddl]_{\sP^\#\zh}
 & & & \sT \ul{\bA^\#}\ar[dr]^{\sT\zh^\#}\ar[ddl]_/-20pt/{\tau_{\ul{\bA^\#}}}
 & \\
 & & \ul{\bA}\ar[rrr]^/-20pt/{\cE_l}\ar[ddl]_/-20pt/{\ul{\zh}}
 & & & \sT M \ar[ddl]_{\tau_M}\\
 \ul{\bA^\#}\ar[rrr]^/-20pt/{id}\ar[dr]^{\ul{{\zh^\#}}}
 & & & \ul{\bA^\#}\ar[dr]^{\ul{{\zh^\#}}} & &  \\
 & M\ar[rrr]^{id}& & & M &
}\label{affgebroid}
\ee
Every such morphism is the composition of (\ref{iso1}) with a
morphisms of double affine bundles \be\xymatrix{
 & \sP\bA^\#\ar[rrr]^{\wt{\zG}} \ar[dr]^{\sP^\#\zh^\#}
 \ar[ddl]_{\sP\zh^\#}
 & & & \sT \ul{\bA^\#}\ar[dr]^{\sT\zh^\#}\ar[ddl]_/-20pt/{\tau_{\ul{\bA^\#}}}
 & \\
 & & \ul{\bA}\ar[rrr]^/-20pt/{\wt{\zG}_l}\ar[ddl]_/-20pt/{\ul{\zh}}
 & & & \sT M \ar[ddl]_{\tau_M}\\
 \ul{\bA^\#}\ar[rrr]^/-20pt/{id}\ar[dr]^{\ul{{\zh^\#}}}
 & & & \ul{\bA^\#}\ar[dr]^{\ul{{\zh^\#}}} & &  \\
 & M\ar[rrr]^{id}& & & M &
}\label{affgebroid1}
\ee
Such morphisms correspond to affine tensors
$\zG=\zG_\cE\in\Sec(\wt{\sT}\bA^\#\ot_{\ul{\bA^\#}}\sT\ul{\bA^\#})$
or, equivalently, to affine-closed biderivation brackets
\be\label{abd}\{\cdot,\cdot\}_{\zG}:\Sec(\bAP(\bA^\#))\ti
C^\infty(\ul{\bA^\#})\ra C^\infty(\ul{\bA^\#}).\ee This obvious
equivalence is induced by the standard identification
$$\wt{\sT}\bA^\#\ot_{\ul{\bA^\#}}\sT\ul{\bA^\#}\simeq
\Hom_{\ul{\bA^\#}}((\wt{\sT}\bA^\#)^\*,\sT\ul{\bA^\#}) \simeq
\Hom_{\ul{\bA^\#}}(\wh{\sP\bA^\#},\sT\ul{\bA^\#})\simeq
\Aff_{\ul{\bA^\#}}(\sP\bA^\#,\sT\ul{\bA^\#}),$$ so
$$\{\zs,f\}_\zG=\la\zG,\bd\zs\ot\xd f\ran,$$ where the affine de
Rham differential $\bd\zs\in\Sec(\sP\bA^\#)$ is regarded also as a
section of the vector hull $\wh{\sP\bA^\#}$.

That (\ref{abd}) is an affine-closed biderivation means that the
bracket is an (affine) derivation with respect to the first
argument, a derivation with respect to the second argument, and it
is affine-closed, i.e. the bracket $\{ \zs,f\}_{\zG}$ of an affine
section $\zs:\ul{\bA^\#}\ra \bA^\#$ and an affine function $f$ on
$\ul{\bA^\#}$ is an affine function on $\ul{\bA^\#}$.

Note that such brackets are just affine-linear parts of certain
affine-closed biderivations
$$\{\cdot,\cdot\}:\Sec(\bAP(\bA^\#))\ti \Sec(\bAP(\bA^\#))\ra
C^\infty(\ul{\bA^\#}).$$ On the level of tensors it means that
$\zG$ can be understood as the projection of a tensor from
$\Sec(\wt{\sT}\bA^\#\ot_{\ul{\bA^\#}}\wt{\sT}\bA^\#)$, i.e. as the
projection of a $\chi_\bA$-invariant and affine 2-contravariant
tensor on $\bA^\#$. Note that in the skew-symmetric case, e.g. in
the case of an aff-Poisson bracket, there is a one-to-one
correspondence between skew biderivations and their affine-linear
parts, since $\{ a,b\}=\{ a,b-a\}^2_\sv$.

According to the identification of sections of $\bA$ with affine
sections of $\bAP(\bA^\#)$(\cite[Theorem 13]{GGU2}) we can derive
out of the bracket $\{\cdot,\cdot\}_\zG$, similarly as it has been
done for aff-Poisson brackets, a general {\it special affgebroid
bracket}
$$[\cdot,\cdot]_\zG:\Sec(\bA)\ti\Sec(\sV(\bA))\ra\Sec(\sV(\bA)).$$
The following theorem (which is completely analogous to
\cite[Theorem 23]{GGU2}, so we skip the proof) explains  in
details what we understand as a special affgebroid bracket.
\begin{theorem}     \label{gebroid}
A special affgebroid structure (\ref{affgebroid}) can be
equivalently defined as an affine-linear bracket
$$[\cdot,\cdot]_\cE:\Sec(\bA)\ti\Sec(\sV(\bA))\ra\Sec(\sV(\bA))$$
which is  special (i.e.
$[a,u+v_\bA]_\cE=[a+v_\bA,u]_\cE=[a,u]_\cE$), together with an
affine bundle morphisms $\cE_l \colon A\rightarrow \sT M$ and a
vector bundle morphism $\cE_r \colon \sV(A)\rightarrow \sT M$(left
and right anchors), such that
                $$ [a,gY]_\cE = g[a,Y]_\cE+(\cE_l\circ a)(g)Y
                                                   $$
                                                   and
                $$ [a+fX,Y]_\cE = (1-f)[a,Y]_\cE+f[a+X,Y]_\cE-(\cE_r\circ Y)(f)X.
                                                   $$
where $a$ is a section of $\bA$ and $X$ is a section of
$\sV(\bA)$. The brackets $[\cdot,\cdot]_\cE$,
$\{\cdot,\cdot\}_\cE$ and the tensor
$\zG=\zG_\cE\in\Sec(\wt{\sT}\bA^\#\ot\sT\ul{\bA^\#})$ are related
by the formula
$$\la\zG_\cE,\bd\zi^\#_a\ot\xd\zi_X\ran= \zi^\#_{[a,X]_\cE}=\{\zi^\#_a,\zi^\#_X\}_\cE, $$
where $\zi^\#_a$ (resp., $\zi^\#_X$) is the corresponding section
of $\bAP(\bA^\#)$ (resp., the corresponding function on
$\ul{\bA^\#}$). The special affgebroid is a special Lie affgebroid
if and only if the tensor $\zG_\cE$ is an aff-Poisson tensor.
\end{theorem}
In local affine coordinates, every  $\cE$ as above is of the form
                \bea\label{cE} &\cE(x^a,y^i,p_b,\zp_j) =(x^a,
\zp_i,  \\
                &\zr^b_0(x)+\sum_k\zr^b_k(x)y^k,\,
c^m_{0j}(x)+\sum_kc^k_{0j}(x)\zp_k+\sum_ic^m_{ij}(x)y^i+\sum_{i,k}c^k_{ij}(x)
y^i\zp_k + \sum_a\zs^a_j(x) p_a),
                                                             \nn\eea
where $i,j,k=1,\dots,m-1$, and $\cE$ corresponds to the affine
2-contravariant tensor $\zG_\cE$ on $\bA^\#$
                \bea \label{tensor} \zG_\cE =& \sum_{i=0,j=1}^{m-1}\left(c^m_{ij}(x)+
               \sum_{k=1}^{m-1} c^k_{ij}(x)\zx_k\right) \partial _{\zx_i}\otimes \partial
_{\zx_j} +\\& \sum_b\left(\sum_{i=0}^{m-1}\zr^b_i(x) \partial
_{\zx_i} \otimes\partial _{x^b} - \sum_{j=1}^{m-1}\zs^b_j(x)
\partial _{x^b} \otimes\partial _{\zx_j}\right).\nn
            \eea
The corresponding affgebroid bracket on
$$\Sec(\bA)\ti\Sec(\sV(\bA))\subset
\Sec(\wh{\bA})\ti\Sec(\wh{\bA})$$ reads
$$[e_0+\sum_{i=1}^mf_ie_i,\sum_{j=1}^mg_je_j]=\sum_{k=1}^{m}\left(
\sum_{i=0,j=1}^{m-1}f_ig_jc_{ij}^k+\sum_a\left(\sum_{i=0}^{m-1}\zr^a_if_i\frac{\pa
g_k}{\pa x^a}-\sum_{j=1}^{m-1}\zs^a_jg_j\frac{\pa f_k}{\pa
x^a}\right)\right)e_k,$$ with the convention that $f_0=1$.

\subsection{Lagrangian and Hamiltonian formalisms}
Combining (\ref{affgebroid}) and (\ref{affgebroid1}) we get the
{\it affine Tulczyjew triple}:
\be\xymatrix{
 & \sP\bA^\#\ar[rrr]^{\wt{\zG}} \ar[dr]^{\sP^\#\zh^\#}
 \ar[ddl]_{\sP\zh^\#}
 & & & \sT \ul{\bA^\#}\ar[dr]^{\sT\zh^\#}\ar[ddl]_/-20pt/{\tau_{\ul{\bA^\#}}}
 &  &  & \sP \bA \ar[ddl]_/-25pt/{\sP^\#\zh} \ar[dr]^{\sP\zh}
\ar[lll]_{\cE}
 & \\
 & & \ul{\bA}\ar[rrr]^/-20pt/{\wt{\zG}_r}\ar[ddl]_/-20pt/{\ul{\zh}}
 & & & \sT M \ar[ddl]_{\tau_M}
& & & \ul{\bA}\ar[lll]_/+20pt/{\cE_l}\ar[ddl]_{\ul{\zh}}
 \\
 \ul{\bA^\#}\ar[rrr]^/-20pt/{id}\ar[dr]^{\ul{{\zh^\#}}}
 & & & \ul{\bA^\#}\ar[dr]^{\ul{{\zh^\#}}} & & & \ul{\bA^\#}
 \ar[dr]^{\ul{\zh^\#}}\ar[lll]_/-20pt/{id}
 & &  \\
 & M\ar[rrr]^{id}& & & M & & & M\ar[lll]_{id}&
}\label{Tt}
\ee
The left-hand side is Hamiltonian, the right-hand side is
Lagrangian, and the 'dynamics' lives in the middle.

Any Lagrangian, that is a section $\cL:\ul{\bA}\ra\bA$ of the
AV-bundle $\bAP(\bA)$, defines a Lagrangian submanifold $\cN=(\bd
\cL)(\ul{\bA})$ in $\sP\bA$, being the image of the affine de Rham
differential $\bd\cL$, i.e. the image of the section $\bd
\cL:\ul{\bA}\ra\sP\bA$. The further image $\cD=\cE(\cN)$ can be
understood as an implicit differential equation on $\ul{\bA^\#}$,
solutions of which are `phase trajectories' of the system. The
Lagrangian defines also smooth maps: $\cL_{eg}:\ul{\bA}\ra
\ul{\bA^\#}$ -- the Legendre mapping associated with $\cL$, and
$\wt{\cL}_{eg}:\ul{\bA}\ra \sT\ul{\bA^\#}$, by
$$\cL_{eg}=\zt_{\ul{\bA^\#}}\circ\cE\circ\bd\cL=\sP^\#\zh\circ\bd\cL$$ and
$\wt{\cL}_{eg}=\cE\circ\bd\cL$. The map $\cL_{eg}$  is {\it de
facto} the vertical derivative of the section $\cL$. The
introduced ingredients produce an implicit differential equation,
this time for curves $\zg:I\ra\ul{\bA}$.

This equation, which will be denoted by $(E_\cL)$,  is represented
by a subset $E_\cL$ of $\sT\ul{\bA}$ being  the inverse image
$$E_\cL=\sT(\wt{\cL}_{eg})^{-1}(\sT^2\ul{\bA^\#})$$ of the subbundle
$\sT^2\ul{\bA^\#}$ of holonomic vectors in $\sT\sT\ul{\bA^\#}$.
The solutions are such paths $\zg:I\ra\ul{\bA}$ that the tangent
prolongation $\sT(\cL_{eg}\circ\zg)$ of $\cL_{eg}\circ\zg$ is
exactly $\wt{\cL}_{eg}\circ\zg$.

In local coordinates $\cL$ is just a function $\cL=\cL(x,y)$ and
$\cD$ has the parametrization by $(x^a,y^i)$ via $\wt{\cL}$ in the
form (cf. (\ref{cE}))
\bea\label{aL}&\wt{\cL}_{eg}(x^a,y^i)=(x^a,
\frac{\pa\cL}{\pa y^i}(x,y), \zr^b_0(x)+\sum_k\zr^b_k(x)y^k,\\
&c^m_{0j}(x)+\sum_ic^m_{ij}(x)y^i+\sum_k\frac{\pa\cL}{\pa
y^k}(x,y)\left(c^k_{0j}(x)+\sum_ic^k_{ij}(x) y^i\right) +
\sum_a\zs^a_j(x)\frac{\partial\cL}{\partial x^a}(x,y)),\nn
\eea
and the equation $(E_\cL)$, for $\zg(t)=(x^a(t),y^i(t))$, is the
sytem of equations
\be\frac{\xd x^a}{\xd t}=\zr^a_0(x)+\sum_k\zr^a_k(x)y^k,
\label{aEL1}\ee
\be \frac{\xd}{\xd t}\left(\frac{\partial\cL}{\partial
y^j}\right)=c^m_{0j}(x)+\sum_ic^m_{ij}(x)y^i+\sum_k\frac{\pa\cL}{\pa
y^k}\left(c^k_{0j}(x)+\sum_ic^k_{ij}(x) y^i\right) +
\sum_a\zs^a_j(x) \frac{\partial\cL}{\partial x^a}.\label{aEL2}\ee
Note that in the particular case when the special affine bundle is
trivial, $\bA=A_0\ti\bI$, and the special affgebroid structure on
$\bA$ comes from the product of a Lie affgebroid structure on
$A_0$ and the trivial Lie algebroid structure in $\bI$, this is in
the full agreement with \cite[(3.14)]{IMPS}, if only one takes
into account that in this case $\zs^a_j=\zr^a_j$ and $c^m_{ij}=0$.
As one can see from (\ref{aEL1}), the solutions are automatically
admissible curves in $\ul{\bA}$, i.e. the velocity $\frac{\xd
}{\xd t}(\ul{\zh}\circ\zg)(t)$ is $\cE_l(\zg(t))$.

Analogously like in the algebroid case, the Hamiltonian formalism
is related to the tensor $\zG_\cE$ (or the aff-Poisson bracket
$\{\cdot,\cdot\}_\cE$ on $\bAP(\bA^\#)$). By the {\it hamiltonian
vector field} associated with a section $\cH$ of $\bAP(\bA^\#)$ we
understand the vector field on $\ul{\bA^\#}$ associated with the
derivation $\{\cH,\cdot\}_\cE$ of $C^\infty(\ul{\bA^\#})$. Thus
the question of the Hamiltonian description of the dynamics $\cD$
in the simplest form is the question if $\cD$ is the image of a
Hamiltonian vector field. Every such a section $\cH$ we call a
{\it Hamiltonian associated with the Lagrangian} $\cL$.

Like in the algebroid case (cf. Theorem \ref{ham}), when dealing
with a hyperregular Lagrangian section, i.e. when $\cL_{eg}$ is a
diffeomorphism, we can find a Hamiltonian associated with the
Lagrangian $\cL$ explicitly. To describe this "affine Legendre
transformation" let us notice that with every section $\cL$ of
$\bAP(\bA)$ we can associate a map $\wh{\cL}:\ul{\bA}\ra\bA^\#$ as
follows. Let us fix $x\in M$ and $a_x\in(\ul{\bA})_x$ and let
$W_{a_x}$ be the maximal affine subspace in $A_x$ that is tangent
to the submanifold $\cL((\ul{\bA})_x)$ at $\cL(a_x)$. There is a
unique affine function $\wh{\cL}_{a_x}$ on $A_x$ which is from
$\bA^\#_x$ (i.e. $\chi_\bA(\wh{\cL}(a_x))=-1$) and which vanishes
on $W_{a_x}$.

\begin{theo}
If the Lagrangian section $\cL$ is hyperregular, then
$\cH=\wh{\cL}\circ\cL_{eg}^{-1}$ is a section of $\bAP(\bA^\#)$
which is a Hamiltonian associated with $\cL$.
\end{theo}

\noindent {\bf Proof.} Let us use local coordinates and the
pairing $\la\cdot,\cdot\ran_{sa}:\bA\ti\bA^\#\ra\R$ as in
(\ref{par}). Since the distinguished direction in $\bA$ is
$-\pa_{y^m}$, it is easy to see that the affine function
$\wh{\cL}_{a_x}$, $a_x=(x,y^1_0,\dots,y^{m-1}_0)$, is
$$\wh{\cL}_{a_x}(x,y^1,\dots,y^{m-1})=
\left(\sum_{i=1}^{m-1}(y^i-y^i_0)\frac{\pa\cL}{\pa
y^i}(a_x)+\cL(a_x)\right)-y^m,$$ which corresponds to the element
$$\left(x,\sum_{i=1}^{m-1}y^i_0\frac{\pa\cL}{\pa
y^i}(a_x)-\cL(a_x),\frac{\pa\cL}{\pa
y^i}(a_x)\right)\in(\bA^\#)_x,$$ so that
$$\wh{\cL}(x,y)=
\left(x,\sum_{i=1}^{m-1}y^i\frac{\pa\cL}{\pa
y^i}(x,y)-\cL(x,y),\frac{\pa\cL}{\pa y^i}(x,y)\right).$$ Since
$(x,y)\mapsto(x,\frac{\pa\cL}{\pa y^i}(x,y))$ is the Legendre map
$\cL_{eg}$, the composition $\cH=\wh{\cL}\circ\cL_{eg}^{-1}$ is a
section of $\bAP(\bA^\#)$, so
$$\cH(x,\zx)=\left(x,\sum_{i=1}^{m-1}y^i(\zx)\zx_i-\cL(x,y(\zx)),\zx\right)$$
and we end up with the standard Legendre transform.

\hfill{$\blacksquare$}

\section{Examples}
    \begin{example} For an AV-bundle
    $\bZ=(Z,v_\Z)$ take as the lagrangian bundle
the AV-bundle $\bAP(\bA)$ over $\sT M$ with the special affine
(this time, in fact, vector) bundle $\bA=\wt{\sT}\Z$. Such
situation we encounter in the analytical mechanics of a
relativistic charged particle (\cite{TU1}) and in the homogeneous
formulation of Newtonian analytical mechanics. We have here
$\bA^\# = \sP \bZ \times \bI$, $\sP\bA^\# = \sT^\*\sP\bZ$, and
$$ \cE\colon \sP \bA = \sP \wt{\sT}\Z \rightarrow \sT\sP\Z. $$
is the canonical isomorphism \cite{Ur1}. Since a choice of a
section of $\Z$ gives a `linearization' $\Z\simeq M\ti\R$, so that
$\sP\Z\simeq\sT^\* M$ and $\sP\wt{\sT}\Z\simeq\sT^\*\sT M$, we get
the canonical Tulczyjew isomorphism $\cE:\sT^\*\sT M\ra\sT\sT^\*
M$ and, in local coordinates, the classical Euler-Lagrange
equation. The point here is that we use the correct geometrical
object $\wt{\sT}\Z$ which does not refer to any {\it ad hoc}
choice of a section of $\Z$.
    \end{example}

\begin{example}
The Newtonian space-time $N$ is a a system $(N,\zt,g)$, where $N$
is a four-dimensional affine space with the model vector space
$V$, together with the time projection $\zt:V\ra\R$ represented by
a non-zero element of $V^\*$, and an Euclidean metric on $V_0 =
\zt^{-1}(0)\subset V$ represented by a linear isomorphism $g\colon
V_0 \rightarrow V_0^\*$.

It is known that that the standard framework for analytical
mechanics is not appropriate for Newtonian analytical mechanics.
It is useful for the frame-dependent formulation of the dynamics
only. Because of the Newtonian relativity principle which states
that the physics is the same for all inertial observers, the
velocity, the momentum, and the kinetic energy have no vector
interpretation, as for example the sum of velocities depends
strongly on the observer (frame). The equivalence of inertial
frames means that all the above concepts are affine in their
nature and that they become vectors only after fixing an inertial
frame. Of course, to get explicit equations for the dynamics we
usually fix a frame, but a correct geometrical model should be
frame-independent.

In \cite{GU} (see also \cite[Example 11]{GGU2}) an affine
framework for a frame-independent formulation of the dynamics has
been proposed. The presented there construction leads to a special
affine space $\bA_1$, for which $\ul{\bA_1} = V_1 = \zt^{-1}(1)$,
and which is equipped with an affine metric, i.e. a mapping
        $$h \colon \ul{\bA_1}\rightarrow \ul{\bA_1^\#}, $$
with the linear part equal to $mg$, where $m$ is the mass of the
particle. The lagrangian bundle is $\bA = N \times \bA_1$. The
kinetic part of a lagrangian is a unique, up to a constant,
section $\ell$ of $\bAP(\bA)$ such that the Legendre mapping
$\ell_{eg}$ equals $h$. Let $P_1 =\ul{\bA_1^\#}$.  We have the
following, obvious equalities:
\begin{itemize}
\item $\sP\mathbf{A}= N\times V_1\times V^\* \times P_1\simeq
N\times P_1\times V^\* \times V_1=\sP\mathbf{A}^\#,$
\item  $\sT \underline{\mathbf{A}^\#} =  N\times P_1 \times V
\times V_0^\*$ .
\end{itemize}
With this identities,  the mappings $\cE\colon \sP\bA \rightarrow
\sT\ul{\bA^\#} $  and $\wt\zG \colon \sP\bA^\# \rightarrow
\sT\ul{\bA^\#}$ read
\beas
\cE &\colon& \sP\mathbf{A}= N\times V_1\times V^\* \times
P_1\ni(x,v,a,p)
\mapsto (x,p,v,\overline{a}) \in  N\times P_1 \times V \times V_0^\* = \sT\ul{\bA^\#}, \\
\wt\zG &\colon& \sP\bA^\#= N\times P_1\times V^\* \times V_1 \ni
(x,p,a,v) \mapsto (x,p,v, - \overline{a})\in  \sT\ul{\bA^\#},
            \eeas
where $\overline{a}$ is the image of $a$ with respect to the
canonical projection $V^\* \rightarrow V^\*_0$. Now, let us
consider a lagrangian of the form
$$\cL \colon N\times E_1\rightarrow A \colon (x,v)\mapsto (x,\ell(v) - \zf(x)), $$
where $\zf$ is a (time-dependent) potential. We have then
$$\bd\cL \colon N\times E_1\rightarrow \sP \bA \colon (x,v) \mapsto (x,v,-\xd
\zf(x),h(v)) \in \sP\bA $$
    and
$$\sT(\cE \circ \bd \cL)\colon (x,v,v', v'')\mapsto (x,h(v),v, -
\overline{\xd \zf(x)},v',mg(v''), v'',-\sT \overline{\xd
\zf}(x,v')).$$ It follows that
$$E_\cL = \{(x,v,v',v'')\colon v=v', \ \ -\overline{\xd \zf(x)}= mg(v'') \},$$
i.e. the Euler-Lagrange equations read
$$\zt(\dot{x})=1,\quad \ddot{x}=-\frac{1}{m}\nabla\zf(x),$$
where $\nabla\zf(x)=g^{-1}(\overline{\xd \zf(x)})$ is the "space"
gradient of $\zf$ at $x\in N$.
    \end{example}


\bigskip
\noindent Katarzyna Grabowska\\
Division of Mathematical Methods in Physics \\
                University of Warsaw \\
                Ho\.za 69, 00-681 Warszawa, Poland \\
                 {\tt konieczn@fuw.edu.pl} \\\\
\noindent Janusz Grabowski\\Polish Academy of Sciences\\Institute
of Mathematics\\\'Sniadeckich 8, P.O. Box 21, 00-956 Warszawa,
Poland\\{\tt jagrab@impan.gov.pl}\\\\
\noindent Pawe\l\ Urba\'nski\\
Division of Mathematical Methods in Physics \\
                University of Warsaw \\
                Ho\.za 69, 00-681 Warszawa, Poland \\
                 {\tt urba\'nski@fuw.edu.pl}

\end{document}